\documentclass[conference,letterpaper]{IEEEtran}
\usepackage{fancyhdr}
\usepackage{accents}
\usepackage{amssymb}
\usepackage{graphicx}
\usepackage{setspace}
\usepackage{mathrsfs}

\newcommand{\bfa}{\mathbf{a}}
\newcommand{\bfb}{\mathbf{b}}

\newcommand{\bff}{\mathbf{f}}
\newcommand{\bfg}{\mathbf{g}}

\newcommand{\bfu}{\mathbf{u}}
\newcommand{\bfv}{\mathbf{v}}
\newcommand{\bfx}{\mathbf{x}}

\newcommand{\bfA}{\mathbf{A}}

\newcommand{\bfF}{\mathbf{F}}

\newcommand{\bfI}{\mathbf{I}}

\newcommand{\bfL}{\mathbf{L}}
\newcommand{\bfM}{\mathbf{M}}

\newcommand{\bfQ}{\mathbf{Q}}
\newcommand{\bfR}{\mathbf{R}}

\newcommand{\bfU}{\mathbf{U}}
\newcommand{\bfV}{\mathbf{V}}

\newcommand{\bfTh}{\mathbf{\Theta}}
\newcommand{\zeros}{\mathbf{0}}
\newcommand{\sM}{\mathcal{M}}

\newcommand{\expect}{\mathbb{E}}
\newcommand{\expectx}{\mathbb{E}_{\bfx}}

\newcommand{\bffg}{\accentset{\frown}{\mathbf{f}}}
\newcommand{\bfxg}{\accentset{\frown}{\mathbf{x}}}
\newcommand{\bfyg}{\accentset{\frown}{\mathbf{y}}}

\usepackage{fancyhdr}
\newtheorem{lemma}{Lemma}
\newtheorem{theorem}{Theorem}
\newtheorem{corollary}{Corollary}

\begin{document}

\title{Analysis of Discrete and Hybrid Stochastic Systems by Nonlinear
  Contraction Theory}

\author{\authorblockN{Quang-Cuong Pham}
\authorblockA{Laboratoire de Physiologie de la Perception et de l'Action\\
Coll\`ege de France - CNRS\\
Paris, France\\
cuong.pham@normalesup.org}}

\maketitle

\thispagestyle{fancy}
\fancyhead{}
\lhead{}
\cfoot{}
\rfoot{}
\renewcommand{\headrulewidth}{0pt}
\renewcommand{\footrulewidth}{0pt}

\begin{abstract}
  We investigate the stability properties of discrete and hybrid
  stochastic nonlinear dynamical systems. More precisely, we extend
  the stochastic contraction theorems (which were formulated for
  continuous systems) to the case of discrete and hybrid resetting
  systems. In particular, we show that the mean square distance
  between any two trajectories of a discrete (or hybrid resetting)
  contracting stochastic system is upper-bounded by a constant after
  exponential transients. Using these results, we study the
  synchronization of noisy nonlinear oscillators coupled by discrete
  noisy interactions.
\end{abstract}


\section{Introduction}

Contraction theory is a set of relatively recent tools that provide a
systematic approach to the stability analysis of a large class of
nonlinear dynamical systems \cite{LS98,WS05,ES06tac,PhaX07}. A
nonlinear nonautonomous system $\dot{\bfx}=\bff(\bfx,t)$ is
\emph{contracting} if the symmetric part of the Jacobian matrix of
$\bff$ is uniformly negative definite in some metric. Using elementary
fluid dynamics techniques, it can be shown that contracting systems
are \emph{incrementally stable}, that is, any two system trajectories
exponentially converge to each other~\cite{LS98}.

From a practical viewpoint, contraction theory has been successfully
applied to a number of important problems, such as mechanical
observers and controllers design~\cite{LS00mech}, chemical processes
control~\cite{LS00chem}, synchronization analysis \cite{WS05,PS07} or
biological systems modelling \cite{GirX08}.

Recently, contraction analysis has been extended to the case of
\emph{stochastic} dynamical systems governed by It\^o differential
equations \cite{PhaX07}. In parallel, hybrid versions of contraction
theory have also been developped \cite{ES06tac}. A hybrid system is
characterized by a \emph{continuous} evolution of the system's state,
and intermittent \emph{discrete} transitions. Such systems are
pervasive in both artificial (e.g. analog physical processes
controlled by digital devices) and natural (e.g. spiking neurons with
subthreshold dynamics) environments.

This paper benefits from these recent developments, and provides an
exponential stability result for discrete and hybrid systems governed
by stochastic \emph{difference} and \emph{differential}
equations. More precisely, we prove in section \ref{sec:discrete} and
\ref{sec:hybrid} that the mean square distance between any two
trajectories of a discrete (respectively hybrid resetting) stochastic
contracting system is upper-bounded by a constant after exponential
transients. This bound can be expressed as function of the noise
intensities and the contraction rates of the noise-free systems. In
section~\ref{sec:remarks}, we briefly discuss a number of theoretical
issues regarding our analysis. In section \ref{sec:example}, we study,
using the previously developped tools, the synchronization of noisy
nonlinear oscillators that interact by discrete noisy
couplings. Finally, some future directions of research are indicated
in section \ref{sec:conclusion}.

\textbf{Notations} The symmetric part of a matrix $\bfA$ is defined as
$\bfA_s=\frac{1}{2}\left(\bfA+\bfA^T\right)$. For a symmetric matrix
$\bfA$, $\lambda_{\min}(\bfA)$ and $\lambda_{\max}(\bfA)$ denote
respectively the smallest and the largest eigenvalue of $\bfA$. A set
of symmetric matrices $\left(\bfA_i\right)_{i\in I}$ is uniformly
positive definite if $ \exists \alpha>0,\ \forall i\in I,\
\lambda_{\min}(\bfA_i)\geq \alpha$. Finally, for a process $\bfx(t)$,
we note $\expectx(\cdot)=\expect(\cdot|\bfx(0)=\bfx)$.

\section{Discrete systems}

\label{sec:discrete}

We first prove a lemma that makes explicit the initial ``discrete
contraction'' proof (see section 5 of \cite{LS98}). Note that a
similar proof for continuous systems can be found in \cite{AR03}.

\begin{lemma}[and definition]
  \label{lemma:ineq}
  Consider two metrics $\bfM_i=\bfTh_i^T\bfTh_i$ defined over
  $\mathbb{R}^{n_i}$ ($i=1,2$) and a smooth function
  $\bff:\mathbb{R}^{n_1} \to\mathbb{R}^{n_2}$. The generalized
  Jacobian of $\bff$ in the metrics $(\bfM_1,\bfM_2)$ is defined by
  \[
  \bfF=\bfTh_2\frac{\partial \bff}{\partial \bfx}\bfTh_1^{-1}
  \]

  Assume now that $\bff$ is \emph{contracting} in the metrics
  $(\bfM_1,\bfM_2)$ with rate $\beta$ $(0<\beta<1)$, i.e.
  \[
  \forall \bfx\in\mathbb{R}^{n_1} 
  \quad \lambda_{\max}(\bfF(\bfx)^T\bfF(\bfx)) \leq \beta
  \]
  
  Then for all $\bfu,\bfv\in\mathbb{R}^n$, one has
  \[
  d_{\bfM_2}(\bff(\bfu),\bff(\bfv))^2 \leq \beta d_{\bfM_1}(\bfu,\bfv)^2
  \]
  where $d_{\bfM}$ denotes the distance associated with the metric
  $\bfM$ (the distance between two points is defined by the infimum of
  the lengths in the metric $\bfM$ of all continuously differentiable
  curves connecting these points).
\end{lemma}

\textbf{Proof} Consider a C$^1$ curve $\gamma :
[0,1]\to \mathbb{R}^{n_1}$ that connects $\bfu$ and $\bfv$
(i.e. $\gamma(0)=\bfu$ and $\gamma(1)=\bfv$).  The $\bfM_1$-length of
such a curve is given by
\[
L_{\bfM_1}(\gamma)=\int_0^1 \sqrt{\left(\frac{\partial \gamma}{\partial
    u}(u)\right)^T \bfM_1\left(\frac{\partial \gamma}{\partial
    u}(u)\right)} du
\]

Since $\bff$ is a smooth function, $\bff(\gamma)$ is also a C$^1$
curve, with
\[
L_{\bfM_2}(\bff(\gamma))=\int_0^1 \sqrt{\left(\frac{\partial
      \bff(\gamma)}{\partial u}(u)\right)^T \bfM_2\left(\frac{\partial
      \bff(\gamma)}{\partial u}(u)\right)} du
\]

The chain rule next implies that
\[
\frac{\partial \bff(\gamma)}{\partial u}(u)=\frac{\partial
  \bff}{\partial \bfx} \frac{\partial \gamma}{\partial u}(u)
\]
which leads to
\begin{equation}
  \label{eq:L}
  \begin{array}{rcl}
    L_{\bfM_2}(\bff(\gamma))
    &=&\int_0^1 \sqrt{ \left(
        \frac{\partial \gamma}{\partial u}^T
        \frac{\partial \bff}{\partial \bfx}^T\bfTh_2^T\bfTh_2
        \frac{\partial \bff}{\partial \bfx} 
        \frac{\partial \gamma}{\partial u} \right)} du \\
    &=&\int_0^1 \sqrt{
      \left(\frac{\partial \gamma}{\partial u}^T\bfTh_1^T\right)
      \bfF^T\bfF
      \left(\bfTh_1\frac{\partial \gamma}{\partial u}\right)} du \\
    &\leq& \int_0^1 \sqrt{\beta \left (\frac{\partial \gamma}{\partial u}^T 
        \bfTh_1^T\bfTh_1
        \frac{\partial \gamma}{\partial u} \right)} du \\
    &=&\sqrt{\beta} L_{\bfM_1}(\gamma)
  \end{array}
\end{equation}

Choose now a sequence of curves $(\gamma_n)_{n\in\mathbb{N}}$ such that
$\lim_{n\to \infty}L_{\bfM_1}(\gamma_n)=d_{\bfM_1}(u,v)$. From
(\ref{eq:L}), one has $\forall n\in \mathbb{N},\
L_{\bfM_2}(\bff(\gamma_n))\leq \sqrt{\beta} L_{\bfM_1}(\gamma_n)$. By
definition of distance, one then has $\forall n\in \mathbb{N},\
d_{\bfM_2}(\bff(u),\bff(v))\leq \sqrt{\beta}
L_{\bfM_1}(\gamma_n)$. Finally, by letting $n$ go to infinity in the
last inequality, one obtains the desired result. $\Box$

\begin{theorem}[Discrete stochastic contraction]
  \label{theo:discrete}

  Consider the stochastic difference equation

  \begin{equation}
    \label{eq:main}
    \bfa_{k+1}=\bff(\bfa_k,k)+\sigma(\bfa_k,k)w_{k+1}
  \end{equation}
  where $\bff$ is a $\mathbb{R}^n \times \mathbb{N}\to\mathbb{R}^n$
  function, $\sigma$ is a $\mathbb{R}^n \times
  \mathbb{N}\to\mathbb{R}^{nd}$ matrix-valued function and
  $\{w_k,k=1,2,\dots\}$ is a sequence of independent $d$-dimensional
  Gaussian noise vectors, with $w_k\sim \mathscr{N}(\zeros,\bfQ_k)$.
  
  Assume that the system verifies the following two hypotheses
  
  \begin{description}
  \item[\textbf{(H1)}] the dynamics $\bff(\bfa,k)$ is contracting in
    the metrics $(\bfM_k,\bfM_{k+1})$, with contraction rate $\beta$
    $(0<\beta<1)$, and the metrics $(\bfM_k)_{k\in\mathbb{N}}$ are
    uniformly positive definite.

  \item[\textbf{(H2)}] the impact of noise is uniformly upper-bounded
    by a constant $\sqrt{C}$ in the metrics $\bfM_k$
    \[
    \forall \bfa,k \quad
    d_{\bfM_k}(\bff(\bfa,k),\bff(\bfa,k)+\sigma(\bfa,k)w_k)\leq
    \sqrt{C}
    \]
 
  \end{description}
  
  Let $\bfa_k$ and $\bfb_k$ be two trajectories whose initial
  conditions are given by a probability distribution
  $p(\bfx_0)=p(\bfa_0,\bfb_0)$. Then for all $k\geq 0$
  \begin{eqnarray}
    \label{eq:etot}
    \expect \left(d_{\bfM_k}(\bfa_k,\bfb_k)\right) \leq
    \frac{2\sqrt{C}}{1-\sqrt{\beta}}  + \nonumber \\
    \sqrt{\beta}^k \int
    \left[d_{\bfM_0}(\bfa,\bfb)-\frac{2\sqrt{C}}{1-\sqrt{\beta}}\right]^+
    dp({\bfa},{\bfb})
  \end{eqnarray}
  where $[\cdot]^+=\max(0,\cdot)$. 

  This implies in particular that for all $k \geq 0$
  \begin{equation}
    \label{eq:e} 
     \expect \left(d_{\bfM_k}(\bfa_k,\bfb_k) \right) \leq
     \frac{2\sqrt{C}}{1-\sqrt{\beta}} + 
     \sqrt{\beta}^k\expect\left(d_{\bfM_0}(\bfa_0,\bfb_0)\right)
  \end{equation}  
\end{theorem}

\textbf{Proof} Let $\bfx=(\bfa,\bfb)^T\in\mathbb{R}^{2n}$. We have by
the triangle inequality (to avoid long formulas, we drop the second
argument of $\bff$ and $\sigma$ in the following calculations)
\[
\begin{array}{rcl}
d_{\bfM_{k+1}}(\bfa_{k+1},\bfb_{k+1})
&\leq& d_{\bfM_{k+1}}(\bff(\bfa_k),\bff(\bfb_k))\\
&+&d_{\bfM_{k+1}}(\bff(\bfa_k),\bff(\bfa_k)+\sigma(\bfa_k) w_{k+1})\\
&+&d_{\bfM_{k+1}}(\bff(\bfb_k),\bff(\bfb_k)+\sigma(\bfb_k) w'_{k+1})
\end{array}
\]

Let us examine the conditional expectations of the three terms of the
right hand side

\begin{itemize}
\item From \textbf{(H1)} and lemma \ref{lemma:ineq} one has
  \[
  \expectx (d_{\bfM_{k+1}}(\bff(\bfa_k),\bff(\bfb_k))) \leq
  \sqrt{\beta} \expectx (d_{\bfM_k}(\bfa_k,\bfb_k))
  \]
\item Next, from \textbf{(H2)}
  \[
  \expectx(d_{\bfM_{k+1}}(\bff(\bfa_k),\bff(\bfa_k)+\sigma(\bfa_k)
  w_{k+1}))  \leq \sqrt{C}
  \]
  and similarly for
  $d_{\bfM_{k+1}}(\bff(\bfb_k),\bff(\bfb_k)+\sigma(\bfb_k)
  w'_{k+1})$. 
\end{itemize}

If we now set
$u_k=\expectx(d_{\bfM_k}(\bfa_k,\bfb_k))$ then
the above implies
\begin{equation}
  \label{eq:u}
  u_{k+1}\leq \sqrt{\beta} u_k+2\sqrt{C}
\end{equation}

Define next $v_k=u_k-2\sqrt{C}/(1-\sqrt{\beta})$. Then replacing $u_k$
by $v_k+2\sqrt{C}/(1-\sqrt{\beta})$ in (\ref{eq:u}) yields
\[
v_{k+1}\leq\sqrt{\beta} v_k
\]

This implies that $\forall k\geq 0, \ v_k\leq v_0\sqrt{\beta}^k \leq
[v_0]^+\sqrt{\beta}^k$. Replacing $v_k$ by its expression in terms of
$u_k$ then yields
\[
\forall k\geq 0\quad u_k\leq
\frac{2\sqrt{C}}{1-\sqrt{\beta}}+\sqrt{\beta}
^k\left[u_0-\frac{2\sqrt{C}}{1-\sqrt{\beta}}\right]^+
\]
which is the desired result.

Next, integrating the last inequality with respect to $\bfx$
leads to (\ref{eq:etot}).  Finally, (\ref{eq:e}) follows
from (\ref{eq:etot}) by remarking that
\begin{eqnarray*}
  \label{eq:trick}
  \int
  \left[d_{\bfM_0}(\bfa,\bfb)-\frac{\sqrt{C}}{1-\sqrt{\beta}}\right]^+
  dp(\bfa,\bfb) \leq \\
  \int
  d_{\bfM_0}(\bfa,\bfb) dp(\bfa,\bfb)
  = \expect\left(d_{\bfM_0}(\bfa,\bfb)\right) \quad \Box
\end{eqnarray*}

\textbf{Remark} In the particular context of state-independent
metrics, hypothesis \textbf{(H2)} is equivalent to the following
simpler condition
\[
\forall \bfa,k \quad
\mathrm{tr}\left(\sigma(\bfa,k)^T\bfM_{k+1}\sigma(\bfa,k)\bfQ_k\right)
\leq C
\]

Also, for state-independent metrics, one has
\[
d_{\bfM_k}(\bfa_k,\bfb_k)^2= 
\|\bfa_k-\bfb_k\|_{\bfM_k}^2=
(\bfa_k-\bfb_k)^T\bfM_k(\bfa_k-\bfb_k)
\]
which leads to the following stronger result instead of (\ref{eq:e})
\[
\expect \left(\|\bfa_k-\bfb_k\|_{\bfM_k}^2 \right) \leq
\frac{2C}{1-\beta} + 
\beta^k\expect\left(\|\bfa_0-\bfb_0\|_{\bfM_0}^2\right)
\]

\section{Hybrid systems}

\label{sec:hybrid}

We have derived above the discrete stochastic contraction theorem for
\emph{time- and state-dependent} metrics, contrary to the context of
continuous systems, where the state-dependent-metrics version of the
contraction theorem is still unproved~\cite{PhaX07}. We now address
the case of hybrid systems, but due to the current limitations of
continuous stochastic contraction, only state-independent metrics will
be considered.

For clarity, we assume in this paper \emph{constant dwell-times},
although more elaborate conditions regarding dwell-times can be
adapted from \cite{ES06tac}.

Consider the hybrid resetting stochastic dynamical system
\begin{equation}
  \label{eq:disc}
  \forall k\geq 0 \quad
  \bfa(k\tau^+)=\bff_d(\bfa(k\tau^-),k)+\sigma_d(\bfa(k\tau^-),k)w_k
\end{equation}
\begin{equation}
  \label{eq:cont}
  \forall k\geq 0,\ \forall t\in ]k\tau,(k+1)\tau[\quad
  d\bfa=\bff_c(\bfa,t)dt+\sigma_c(\bfa,t)dW
\end{equation}

All the contraction properties below will be stated with respect to a
uniformly positive definite time-varying metric
$\bfM(t)=\bfTh(t)^T\bfTh(t)$. Furthermore, it will be assumed that for
all $k\geq 0$, $\bfM$ is continuously differentiable in
$]k\tau,(k+1)\tau[$. Finally, $\bfM(k\tau^-)$ and $\bfM(k\tau^+)$ will
respectively denote the left and right limits of $\bfM(t)$ at
$t=k\tau$ (and similarly for $\bfTh$).

\subsection{The discrete and continuous parts are both contracting}

\begin{theorem}[Hybrid stochastic contraction]
  \label{theo:hybrid-contracting}
  Assume the following conditions
  \begin{description}
  \item[(i)] For all $k$, the discrete part is stochastically contracting
    at $k\tau$ with rate $\beta<1$ and bound ${C_d}$, i.e.
    \[
    \forall \bfa \in\mathbb{R}^n \quad
    \lambda_{\max}\left(\bfF(k\tau)^T\bfF(k\tau)\right) \leq \beta
    \]
    where $\bfF(k\tau)=\bfTh(k\tau^+)
    \frac{\partial \bff_d}{\partial \bfa}(\bfa,k)
    \bfTh(k\tau^-)$, and
    \[
    \forall \bfa\in\mathbb{R}^n \quad
    \mathrm{tr}\left(\sigma_d(\bfa,k)^T\bfM(k\tau^+)\sigma_d(\bfa,k)\bfQ_k\right)
    \leq {C_d}
    \]
  \item[(ii)] For all $k$, the continuous part is stochastically contracting
    in $]k\tau,(k+1)\tau[$ with rate $\lambda>0$ and bound $C_c$, i.e.
    $\forall \bfa\in\mathbb{R}^n,\ \forall t\in]k\tau,(k+1)\tau[,$
    \begin{equation}
      \label{eq:cond-cont}
      \lambda_{\max} \left(\left(\frac{d}{dt}\bfTh(t)+\bfTh(t)\frac{\partial
            \bff}{\partial \bfa}\right)\bfTh^{-1}(t)\right)_s \leq -\lambda
    \end{equation}
    \[
    \mathrm{tr}\left(\sigma_c(\bfa,t)^T\bfM(t)\sigma_c(\bfa,t)\right)
    \leq {C_c}
    \]
  \end{description}

  Let $\bfa(t)$ and $\bfb(t)$ be two trajectories whose initial
  conditions are given by a probability distribution
  $p(\bfx(0))=p(\bfa(0),\bfb(0))$. Then for all $t\geq 0$
  \[
  \begin{array}{c}
    \expect \left(\|\bfa(t)-\bfb(t)\|_{\bfM(t)}^2 \right)  \leq \\
    C_1+
    \expect \left(\|\bfa(0)-\bfb(0)\|_{\bfM(0)}^2 \right)
    \beta^{\lfloor t/\tau \rfloor}e^{-2\lambda t}
  \end{array}
  \]
  where $C_1=\frac{2\lambda C_d+(1-\beta)(1+\beta-r_1)C_c}
  {\lambda(1-\beta)(1-r_1)}$ and $r_1=\beta e^{-2\lambda\tau}$.
\end{theorem}

\textbf{Proof} For all $t\geq 0$, let
$u(t)=\expect\left(\|\bfa(t)-\bfb(t)\|_{\bfM(t)}^2\right)$ and let us
study the evolution of $u(t)$ between $k\tau^+$ and $(k+1)\tau^+$.

Condition (ii) and theorem 2 of \cite{PhaX07} yield
\begin{equation}
  \label{ineq:cont}
  \begin{array}{c}
  u((k+1)\tau^-) \leq \frac{{C_c}}{\lambda} +
  u(k\tau^+) e^{-2\lambda\tau}    
  \end{array}
\end{equation}

Next, condition (i) and theorem~\ref{theo:discrete} above yield
\begin{equation}
  \label{ineq:disc}
  \begin{array}{c}
    u((k+1)\tau^+) \leq \frac{2{C_d}}{1-\beta} +
    \beta u((k+1)\tau^-)
  \end{array}
\end{equation}

Substituting (\ref{ineq:cont}) into (\ref{ineq:disc}) leads to
\[
  \begin{array}{rcl}
    u((k+1)\tau^+)   &\leq& 
    \frac{2{C_d}}{1-\beta} +
    \beta\left(  \frac{{C_c}}{\lambda} +
      \beta u(k\tau^+)  e^{-2\lambda\tau} \right)\\
      &=&
    \frac{2{C_d}}{1-\beta}+\frac{\beta{C_c}}{\lambda}+
    \beta e^{-2\lambda\tau} u(k\tau^+) \nonumber
  \end{array}
\]

Define $D_1=\frac{2{C_d}}{1-\beta}+\frac{\beta{C_c}}{\lambda}$ and
$v_k=u(k\tau^+)-D_1/(1-r_1)$.  Then, similarly to the proof of
theorem~\ref{theo:discrete}, we have $v_{k+1}\leq r_1v_k$, and then
$v_k\leq r_1^k [v_0]^+$, which implies
\begin{eqnarray}
  \label{eq:2}
  u(k\tau^+)&\leq& \frac{D_1}{1-r_1} + \left[u(0^+)-\frac{D_1}{1-r_1}\right]^+ 
  r_1^k \nonumber\\
  &\leq&\frac{D_1}{1-r_1}+u(0^+) r_1^k  \nonumber
\end{eqnarray}

Now, for any $t\geq 0$, choose $k=\lfloor t/\tau \rfloor$. Then
\begin{eqnarray}
  \label{eq:1}
  u(t) &\leq&
  \frac{{C_c}}{\lambda} +
  u(k\tau^+) e^{-2\lambda (t-k\tau)} \nonumber\\
  &\leq& \frac{{C_c}}{\lambda} + \frac{D_1e^{-2\lambda (t-k\tau)}}{1-r_1}
  +  u(0^+) \beta^k e^{-2\lambda t }  \nonumber \\
  &\leq& \frac{{C_c}}{\lambda} + \frac{D_1}{1-r_1}+
  u(0^+) \beta^k e^{-2\lambda t } \nonumber
\end{eqnarray}
which leads to the desired result after some algebraic
manipulations. $\Box$

\subsection{Only the discrete part is contracting}

Let us examine now the more interesting case when the continuous part
is not contracting, more precisely when \mbox{$\lambda\leq 0$} in
(\ref{eq:cond-cont}). For this, we shall need to revisit
the proof of theorem~2 in \cite{PhaX07}.

\begin{theorem}[Case $\lambda=0$]
  \label{theo:hybrid-indiff}
  Assume all the hypotheses of theorem~\ref{theo:hybrid-contracting}
  except that $\lambda=0$ in (\ref{eq:cond-cont}). Then for
  all $t\geq 0$
  \[
  \begin{array}{c}
    \expect \left(\|\bfa(t)-\bfb(t)\|_{\bfM(t)}^2 \right)  \leq \\
    C_2 + \expect \left(\|\bfa(0)-\bfb(0)\|_{\bfM(0)}^2 \right)
    \beta^{\lfloor t/\tau \rfloor}
  \end{array}
  \]
  where $C_2= \frac{2C_d+2\beta(1-\beta)C_c\tau}{(1-\beta)^2}$.
\end{theorem}

\textbf{Proof} As in the proof of theorem~2 in \cite{PhaX07}, let
\[
V(\bfx,t)=V((\bfa,\bfb)^T,t)=(\bfa-\bfb)^T\bfM(t)(\bfa-\bfb)
\]
Lemma 1 of \cite{PhaX07} is unchanged, yielding (see \cite{PhaX07} for
more details)
\[
\forall t\in ]k\tau,(k+1)\tau[ \quad \widetilde{A}V(\bfx(t),t) \leq
2C_c 
\]
where $\widetilde{A}$ is the infinitesimal operator associated with
the process $\bfx(t)$ (see section 2.1.2 of \cite{PhaX07} or p. 15 of
\cite{Kus67} for more details).

By Dynkin's formula \cite{Kus67}, one then obtains for all
$\bfx\in\mathbb{R}^{2n}$
\[
\begin{array}{rcl}
  \expectx V(\bfx(t),t)-V(\bfx,k\tau^+)&=&\expectx
  \int_{k\tau}^t\widetilde{A}V(\bfx(s),s)ds \nonumber \\
  &\leq&\expectx \int_{k\tau}^t 2C_c ds \nonumber \\
  &=& 2C_c (t-k\tau)
\end{array}
\]
Integrating the above inequality with respect to $\bfx$ then yields
\[
\forall t\in ]k\tau,(k+1)\tau[ \quad u(t) \leq 2C_c(t-k\tau) + u(k\tau^+)
\]
In particular, (\ref{ineq:cont}) becomes
\[
u((k+1)\tau^-) \leq 2C_c\tau + u(k\tau^+)
\]
which leads to, after substition into (\ref{ineq:disc}),
\[
u((k+1)\tau^+) \leq \frac{2{C_d}}{1-\beta} + 2\beta C_c\tau +\beta
u(k\tau^+) 
\]

This finally implies 
\[
u(k\tau^+) \leq \frac{\frac{2{C_d}}{1-\beta} + 2\beta
  C_c\tau}{1-\beta} + u(0^+) \beta^k
\]
The remainder of the proof can be adapted from that of
theorem~\ref{theo:hybrid-contracting}. $\Box$

\begin{theorem}[Case $\lambda<0$]
  \label{theo:hybrid-non-contracting}
  Assume all the hypotheses of theorem~\ref{theo:hybrid-contracting}
  except that $\lambda<0$ in (\ref{eq:cond-cont}). Let
  $k=\lfloor t/\tau\rfloor$. There are two cases:
  \begin{itemize}
  \item If $\beta< e^{-2|\lambda|\tau}$, then let $r_2=\beta
    e^{2|\lambda|\tau}<1$. For all $t\geq 0$
    \[
    \begin{array}{c}
      \expect \left(\|\bfa(t)-\bfb(t)\|_{\bfM(t)}^2 \right)  \leq \\
       C_3+\expect \left(\|\bfa(0)-\bfb(0)\|_{\bfM(0)}^2 \right)
      e^{2|\lambda|\tau}r_2^k
    \end{array}
    \]    
    where $C_3=\frac{2|\lambda| C_d+
      (1-\beta)(1+\beta-r_2)e^{2|\lambda|\tau}C_c}
    {|\lambda|(1-\beta)(1-r_2)}$.
  \item If $\beta\geq e^{-2|\lambda|\tau}$, then there is -- in
    general -- no finite bound on $\expect
    \left(\|\bfa(t)-\bfb(t)\|_{\bfM(t)}^2 \right)$ as $t\to +\infty$.
  \end{itemize}
\end{theorem}

\textbf{Proof} One has now for all $t\in ]k\tau,(k+1)\tau[$,
\[
\widetilde{A}V(\bfx(t),t) \leq 2|\lambda| V(\bfx(t),t)+ 2C_c
\]
with $|\lambda|>0$. By Dynkin's formula, one has, for all
$\bfx\in\mathbb{R}^{2n}$
\[
\expectx V(\bfx(t),t)-V(\bfx,k\tau^+) \leq
\expectx\int_{k\tau}^t(2|\lambda| V(\bfx(s),s)+ 2C_c)ds
\]
Let now $g(t)=\expectx V(\bfx(t),t)$. The above equation then yields
\[
g(t)=V(\bfx,k\tau^+)+
2C_c(t-k\tau)+2|\lambda|\int_{k\tau}^t g(s)ds
\]

Applying the classical Gronwall's lemma \cite{Rob06} to $g(t)$ leads
to
\[
\begin{array}{rcl}
  g(t)&\leq& V(\bfx,k\tau^+)+ 2C_c(t-k\tau) +\\
  &&2|\lambda|\int_{k\tau}^t \left(V(\bfx,k\tau^+)+2C_cs\right) \exp\left(\int_s^t
    2|\lambda| du\right)ds\\
  &=&\frac{C_c}{|\lambda|}\left(e^{2|\lambda|(t-k\tau)}-1\right)+
  V(\bfx,k\tau^+)e^{2|\lambda|(t-k\tau)}
\end{array}
\]
Integrating the above inequality with respect to $\bfx$ then yields
$\forall t\in ]k\tau,(k+1)\tau[$,
\[
u(t) \leq
\frac{C_c}{|\lambda|}\left(e^{2|\lambda|(t-k\tau)}-1\right)+
u(k\tau^+)e^{2|\lambda|(t-k\tau)}
\]
which implies
\begin{equation}
  \label{eq:non-contracting}
  u((k+1)\tau^+) \leq D_2
  +\beta e^{2|\lambda|\tau}u(k\tau^+)
  \ 
\end{equation}
where $D_2=\frac{2{C_d}}{1-\beta}+
\frac{\beta C_c}{|\lambda|}\left(e^{2|\lambda|\tau}-1\right)$. 

There are three cases:
\begin{itemize}
\item If $\beta< e^{-2|\lambda|\tau}$, then $r_2=\beta
  e^{2|\lambda|\tau}<1$. By the same reasoning as in
  theorem~\ref{theo:discrete}, one obtains
  \[
  u(k\tau^+) \leq \frac{D_2}{1-r_2} + u(0^+)r_2^k
  \] 
  The remainder of the proof can be adapted from that of
  theorem~\ref{theo:hybrid-contracting}
\item If $\beta=e^{-2|\lambda|\tau}$, then (\ref{eq:non-contracting})
  reads
  \[
  u((k+1)\tau^+) \leq D_2 + u(k\tau^+)
  \]
  which implies $\forall k\geq 0,\ u(k\tau^+)\leq kD_2+u(0^+)$. From
  this, it is clear that there is -- in general -- no finite bound for
  $u(k\tau^+)$.
\item If $\beta> e^{-2|\lambda|\tau}$, then $r_2=\beta
  e^{2|\lambda|\tau}>1$. By the same reasoning as in
  theorem~\ref{theo:discrete}, one obtains
  \[
  u(k\tau^+) \leq \left(u(0^+)+\frac{D_2}{r_2-1}\right) r_2^k
  -\frac{D_2}{r_2-1}
  \]
  Since $r_2>1$ in this case, it is clear that there is -- in general
  -- no finite bound for $u(k\tau^+)$. $\Box$
\end{itemize}

\textbf{Remarks} Theorems \ref{theo:hybrid-indiff} and
\ref{theo:hybrid-non-contracting} show that it is possible to
stabilize an unstable system by discrete resettings. If the continuous
system is \emph{indifferent} ($\lambda=0$), then \emph{any} sequence
of uniformly contracting resettings is stabilizing. However, it should
be noted that the asymptotic bound $C_2\to \infty$ when $\beta\to
1$. In contrast, if the continuous system is \emph{strictly unstable}
($\lambda<0$), then specific contraction rates (depending on the
dwell-time and the ``expansion'' rate of the continuous system) of the
resettings are required. Finally, note that in both cases, the
asymptotic bounds $C_2$ and $C_3$ are increasing functions of the
dwell-time~$\tau$.

\section{Comments}

\label{sec:remarks}

\subsection{Modelling issue: distinct driving noise}
\label{sec:modelling}

In the same spirit as \cite{PhaX07}, and contrary to previous works on
the stability of stochastic systems \cite{Flo95}, the $\bfa$ and
$\bfb$ systems considered in sections \ref{sec:discrete} and
\ref{sec:hybrid} are driven by \emph{distinct} and independent noise
processes. This approach enables us to study the stability of the
system with respect to variations in initial conditions \emph{and} to
random perturbations: indeed, two trajectories of any real-life system
are typically affected by distinct \emph{realizations} of the
noise. In addition, this approach leads very naturally to nice results
regarding the comparison of noisy and noise-free trajectories (see
section \ref{sec:noise-free}), which are particularly useful in
applications (see e.g. section \ref{sec:example}).

However, because of the very fact that the two trajectories are driven
by distinct noise processes, we cannot expect the influence of noise
to vanish when the two trajectories get very close to each other. As a
consequence, the asymptotic bounds $2C/(1-\beta)$ (for discrete
systems) and $C_1$, $C_2$, $C_3$ (for hybrid systems) are strictly
positive. These bounds are nevertheless \emph{optimal}, in the sense
that they can be attained (adapt the Ornstein-Uhlenbeck example in
section 2.3.1 of~\cite{PhaX07}).

\subsection{Noisy and noise-free trajectories}
\label{sec:noise-free}

Instead of considering two noisy trajectories $\bfa$ and $\bfb$ as in
theorem \ref{theo:discrete}, we assume now that $\bfa$ is noisy, while
$\bfb$ is noise-free. More precisely, for all $k\in\mathbb{N}$
\[
\bfa_{k+1}=\bff(\bfa_k,k)+\sigma(\bfa_k,k)w_{k+1}
\]
\[
\bfb_{k+1}=\bff(\bfb_k,k)
\]

To show the exponential convergence of $\bfa$ and $\bfb$ to each
other, one can follow the same reasoning as in the proof of theorem
\ref{theo:discrete}, with $C$ is replaced by $C/2$. This leads to the
following result
\begin{corollary}
  Assume all the hypothesis of theorem \ref{theo:discrete} and
  consider a noise-free trajectory $\bfb_k$ and a noisy trajectory
  $\bfa_k$ whose initial conditions are given by a probability
  distribution $p(\bfa_0)$. Then, for all $k\in\mathbb{N}$
  \begin{eqnarray}
    \label{eq:cor}
    \expect \left(\|\bfa_k-\bfb_k\|_{\bfM_k}^2 \right) \leq
    \frac{C}{1-\beta}  + \nonumber \\
    \beta^k \int
    \left[\|{\bfa}-{\bfb_0}\|_{\bfM_0}^2-\frac{C}{1-\beta}\right]^+
    dp({\bfa})
  \end{eqnarray}
\end{corollary}

\textbf{Remarks}
\begin{itemize}
\item The above derivation of corollary 1 is only permitted by our
  choice of considering distinct driving noise processes for systems
  $\bfa$ and $\bfb$ (see section \ref{sec:modelling}).
\item Based on theorems \ref{theo:hybrid-contracting},
  \ref{theo:hybrid-indiff} and \ref{theo:hybrid-non-contracting},
  similar corollaries can be obtained for hybrid systems.
\item These corollaries provide a robustness result for contracting
  discrete and hybrid systems, in the sense that any contracting
  system is \emph{automatically} protected against noise, as
  quantified by (\ref{eq:cor}). This robustness could be related to
  the exponential nature of contraction stability.
\end{itemize}

\section{Application: oscillator synchronization by discrete
  couplings}

\label{sec:example}

Using the above developped tools, we study in this section the
synchronization of nonlinear oscillators in presence of random
perturbations. The novelty here is that the interactions between the
oscillators occur at \emph{discrete} time instants, contrary to many
previous works devoted to synchronization in the
\emph{state-space}\footnote{Discrete couplings are more frequent in
  the literature devoted to \emph{phase oscillators} synchronization,
  where \emph{phase reduction} techniques are used
  \cite{Izh99}. However, contrary to our approach, these techniques
  are only applicable in the case of weak coupling strenghs and small
  noise intensities.}  \cite{PogX02,PS07}.

Specifically, consider the Central Pattern Generator (CPG) delivering
$2\pi/3$-phase-locked signals of section 5.3 in \cite{PS07}. This CPG
consists of a network of three Andronov-Hopf oscillators
$\bfx_i=(x_i,y_i)^T,\ i=1,2,3$.  We construct below a
discrete-couplings version of this CPG. 

At instants $t=k\tau,\ k\in\mathbb{N}$, the three oscillators are
coupled in the following way (assuming noisy measurements)
\[
\begin{array}{rcl}
 \bfx_i(k\tau^+)&=&\bfx_i(k\tau^-)\\
 &+&\gamma\left(\bfR\left(\bfx_{i+1}(k\tau^-)+\frac{\sigma_d}{\sqrt 2}
    w_k\right) -\bfx_i(k\tau^-)\right)
\end{array}
\]
with $\bfx_4=\bfx_1$ and 
\[
\bfR= \left(\begin{array}{ll}
    -\frac{1}{2} & -\frac{\sqrt{3}}{2}\\
    \frac{\sqrt{3}}{2} & -\frac{1}{2}\\
  \end{array}\right)
\]

Between two interaction instants, the oscillators follow the
uncoupled, noisy, dynamics 
\[
d\bfx_i=\bff(\bfx_i)dt+\frac{\sigma_c}{\sqrt 2} dW
\]
where
\[
\bff(\bfx_i)=\bff\left(\begin{array}{l}x_i\\y_i\end{array}\right)=
\left(\begin{array}{l}
    x_i-y_i-x_i^3-x_iy_i^2\\
    x_i+y_i-y_i^3-y_ix_i^2
\end{array}\right) 
\]

We apply now the projection technique developped in
\cite{PS07,PhaX07}. We recommend the reader to refer to these papers
for more details about the following calculations.

Consider first the (linear) subspace $\sM$ of the global state space
(the global state is defined by $\bfxg=(\bfx_1,\bfx_2,\bfx_3)^T$)
where the oscillators are $2\pi/3$-phase-locked
\[
\sM= \left\{\left(\bfR^2(\bfx),\bfR(\bfx),\bfx\right)^T:
  \bfx\in\mathbb{R}^2\right\}
\]

Let $\bfV$ and $\bfU$ be two orthonormal projections on $\sM^\perp$
and $\sM$ respectively and consider $\bfyg=\bfV\bfxg$. Since the
mapping is linear, using It\^o differentiation rule yields the
following dynamics for $\bfyg$
\begin{equation}
  \label{eq:aux-disc}
  \forall k\in\mathbb{N}\quad
  \bfyg(k\tau^+)=\bfg_d(\bfyg(k\tau^-))+\gamma\frac{\sigma_d}{\sqrt 2} w_k
\end{equation}
\begin{equation}
  \label{eq:aux-cont}
  \forall t\in ]k\tau,(k+1)\tau[\quad
  d\bfyg=\bfg_c(\bfyg)dt+\frac{\sigma_c}{\sqrt 2} dW
\end{equation}
with
\[
\bfg_d(\bfyg)=\bfV\bfL\bfxg=\bfV\bfL(\bfV^T\bfyg+\bfU^T\bfU\bfxg)=\bfV\bfL\bfV^T\bfyg
\]
\[
\bfg_c(\bfyg)=\bfV\bffg(\bfV^T\bfyg+\bfU^T\bfU\bfxg)
\]
where
\[
\bfL=\left(\begin{array}{ccc}
    (1-\gamma)\bfI_2&\gamma\bfR&\zeros \\
    \zeros& (1-\gamma)\bfI_2&\gamma\bfR \\
    \gamma\bfR& \zeros & (1-\gamma)\bfI_2
\end{array}\right)
\]
\[
\bffg(\bfxg)=(\bff(\bfx_1),\bff(\bfx_2),\bff(\bfx_3))^T
\]

Remark that $\bfg_d(\zeros)=\zeros$ and $\bfg_c(\zeros)=\zeros$ (the
last equality holds because of the symmetry of $\bff$: $\forall \bfx,\
\bff(\bfR\bfx)=\bfR(\bff(\bfx))$). Thus, $\zeros$ is a particular
solution to the noise-free version of the hybrid stochastic system
(\ref{eq:aux-disc},\ref{eq:aux-cont}).

Let us now examine the contraction properties of equations
(\ref{eq:aux-disc}) and (\ref{eq:aux-cont}).

We have first
\[
\frac{\partial\bfg_d}{\partial\bfyg}^T\frac{\partial\bfg_d}{\partial\bfyg} 
=\bfV\bfL^T\bfV^T\bfV\bfL\bfV^T=(3\gamma^2-3\gamma+1)\bfI_4
\]
so that $\lambda_{\max}
\left(\frac{\partial\bfg_d}{\partial\bfyg}^T\frac{\partial\bfg_d}{\partial\bfyg}\right)
=3\gamma^2-3\gamma+1<1$ (for $0<\gamma<1$). 

Second, 
\[
\frac{\partial\bfg_c}{\partial\bfyg}=
\bfV\frac{\partial\bffg}{\partial\bfxg}\bfV^T=
\bfV\left(\begin{array}{ccc}
    \frac{\partial\bff}{\partial\bfx}(\bfx_1)&\zeros&\zeros \\
    \zeros& \frac{\partial\bff}{\partial\bfx}(\bfx_2)&\zeros\\
    \zeros& \zeros & \frac{\partial\bff}{\partial\bfx}(\bfx_3)
\end{array}\right)\bfV^T
\]

Now observe that
$\lambda_{\max}\left(\frac{\partial\bff}{\partial\bfx}\right)_s
=1-x^2-y^2\leq 1$. Since $\bfV$ is an orthonormal projection, one then
has
$\lambda_{\max}\left(\frac{\partial\bfg_c}{\partial\bfyg}\right)_s\leq
1$.

Therefore, if
\begin{equation}
  \label{eq:sync-cond}
  3\gamma^2-3\gamma+1<e^{-2\tau}  
\end{equation}
then theorem \ref{theo:hybrid-non-contracting} together with the
corollaries of section \ref{sec:noise-free} imply that, after
exponential transients,
\[
\expect\left(\|\bfyg\|^2\right)\leq 
\frac{2 \gamma^2\sigma_d^2+
  (1-\beta)(1+\beta-\beta e^{2\tau})e^{2\tau}\sigma_c^2}
{2(1-\beta)(1-\beta e^{2\tau})}
\]
where $\beta=3\gamma^2-3\gamma+1$.

To conclude, observe that
\[
\|\bfyg\|^2=\|\bfV\bfxg\|^2=\frac{1}{3}\sum_{i=1}^3\|\bfR\bfx_{i+1}-\bfx_i\|^2
\]

Define the \emph{phase-locking quality} $\delta$ by
\[
\delta= \sum_{i=1}^3\|\bfR\bfx_{i+1}-\bfx_i\|^2
\]
then one finally obtains
\begin{equation}
  \label{eq:final-bound}
  \expect(\delta)\leq 
  \frac{6 \gamma^2\sigma_d^2+3
    (1-\beta)(1+\beta-\beta e^{2\tau})e^{2\tau}\sigma_c^2}
  {2(1-\beta)(1-\beta e^{2\tau})}
\end{equation}
after exponential transients. 

A numerical simulation is provided in Fig. 1.

\begin{figure}[ht]
  \centering
  \includegraphics[scale=0.9]{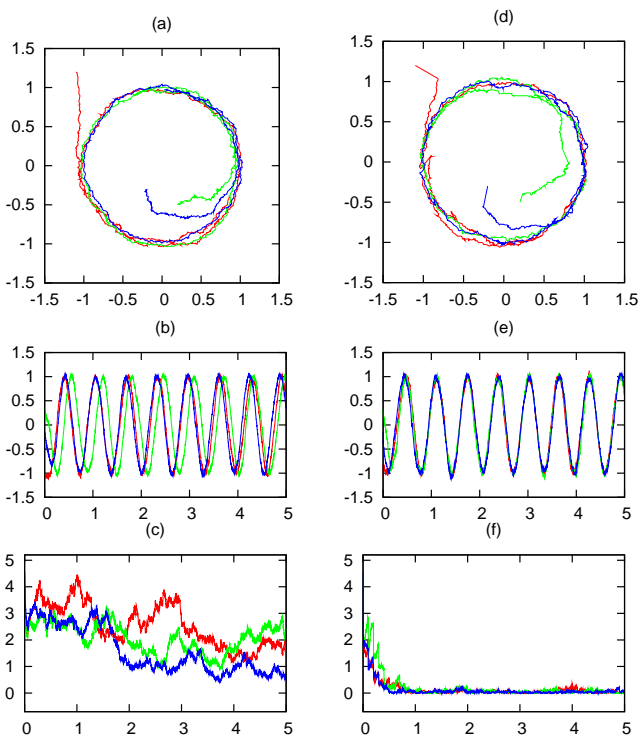}
  \caption{Numerical simulation using the Euler-Maruyama algorithm
    \cite{Hig01}. The following set of parameters was used:
    $\sigma_c=0.1$, $\sigma_d=0.05$, $\tau=0.1$. Two coupling
    strengths were tested: $\gamma_\mathrm{weak}=0.01$ for plots (a),
    (b), (c), and $\gamma_\mathrm{strong}=0.2$ for plots (c), (d),
    (e). Note that $\gamma_\mathrm{weak}$ does not satisfy condition
    (\ref{eq:sync-cond}), while $\gamma_\mathrm{strong}$ does, and
    yields the theoretical bound $\simeq$ 0.446 (as provided by
    (\ref{eq:final-bound})) on the phase-locking quality
    $\delta$. Plots (a) and (d) show the 2d trace of sample
    trajectories of the three oscillators for $t\in[0,1]$. Plots (b)
    and (e) show sample trajectories of the first coordinates of
    $\bfx_1$, $\bfR(\bfx_2)$ and $\bfR^2(\bfx_3)$ as functions of
    time. Plot (c) and (f) show three sample trajectories of~$\delta$.
  }
  \label{fig:simu}
\end{figure}

\section{Perspectives}

\label{sec:conclusion}

We are now focusing on the following directions of research:
\begin{itemize}
\item proving the state-dependent-metrics version of the continuous
  and hybrid stochastic contraction theorems,
\item developping more elaborate conditions on dwell-times, and also
  hybrid \emph{switched} versions of the theorems,
\item applying the synchronization-by-discrete-couplings analysis to
  other types of coupled dynamical systems,
\item studying the robustness of hybrid controllers and observers
  against random perturbations (for instance, the discrete observer
  for inertial navigation developped in~\cite{ZS05}).
\end{itemize}

\section*{Acknowledgment}

The author is grateful to Prof J.-J. Slotine and N. Tabareau for
stimulating discussions, and to Dr H. Hicheur for the careful reading
of the manuscript. This work has been supported by EC - contract
number FP6-IST-027140, action line: Cognitive Systems. This
publication reflects only the author's views. The European Community
is not liable for any use that may be made of the information
contained therein.

\bibliographystyle{IEEEtran}
\bibliography{nsl}

\end{document}